\title{Geography of minimal surfaces of general type with $\mathbb{Z}_2^2$-actions and the locus of Gorenstein stable surfaces.}

\author{Vicente Lorenzo}
\date{}
\RequirePackage{fix-cm}
\RequirePackage{amsmath}

\documentclass{article}
\usepackage{graphicx}
\usepackage[T1]{fontenc}
\usepackage[utf8]{inputenc}
\usepackage{lmodern}
\usepackage[english]{babel}
\usepackage{amssymb}
\usepackage{amsthm}
\usepackage[all]{xy}
\usepackage{enumerate}
\usepackage[none]{hyphenat}
\usepackage{etoolbox}
\usepackage{microtype}
\usepackage{lipsum}
\usepackage[pdfencoding=auto, psdextra]{hyperref}
\usepackage{tikz-cd}

\newtheorem{theorem}{Theorem}

\newtheorem{proposition}{Proposition}

\theoremstyle{remark}\newtheorem{remark}{Remark}
\theoremstyle{remark}\newtheorem{example}{Example}

\newenvironment{acknowledgements}{\textit{Acknowledgements.}}{}

\newcommand\blfootnote[1]{%
  \begingroup
  \renewcommand\thefootnote{}\footnote{#1}%
  \addtocounter{footnote}{-1}%
  \endgroup
}

\begin{document}

\maketitle

\begin{abstract}

In this note the geography of minimal surfaces of general type admitting $\mathbb{Z}_2^2$-actions is studied. More precisely, it is shown that Gieseker's moduli space $\mathfrak{M}_{K^2,\chi}$ contains surfaces admitting a $\mathbb{Z}_2^2$-action for every admissible pair $(K^2, \chi)$ such that $2\chi-6\leq K^2\leq 8\chi-8$ or $K^2=8\chi$.

The examples considered allow to prove that the locus of Gorenstein stable surfaces is not closed in the KSBA-compactification $\overline{\mathfrak{M}}_{K^2,\chi}$ of Gieseker's moduli space $\mathfrak{M}_{K^2,\chi}$ for every admissible pair $(K^2, \chi)$ such that $2\chi-6\leq K^2\leq 8\chi-8$.
\end{abstract}

\section{Introduction.}
Let $X$ be an algebraic surface over the complex numbers, which will be the ground field throughout the paper. The main numerical invariants of $X$ are the self-intersection of its canonical class $K^2_X$ and its holomorphic Euler characteristic $\chi(\mathcal{O}_X)$. 
 If $X$ is minimal and of general type it is well known 
that the following inequalities are satisfied
\begin{equation}\label{GeneralType}
 \chi(\mathcal{O}_X)\geq 1,\quad K_X^2\geq 1, \quad 2\chi(\mathcal{O}_X)-6\leq K_X^2\leq 9\chi(\mathcal{O}_X).
\end{equation}
A pair of integers $(K^2, \chi)$ is said to be an admissible pair if it satisfies the inequalities (\ref{GeneralType}). Whether there exists a minimal algebraic surface of general type with 
given admissible invariants is called the geographical
question and has received the attention of many renowned algebraic geometers. An account of the topic can be found in \cite[Section VII.8.A]{Barth2004}. Several authors have also studied the geography of minimal surfaces of general type with special features. In this note we are interested
in the geography of minimal surfaces of general type admitting a $\mathbb{Z}_2^2$-action. 

Given an admissible pair $(K^2,\chi)$ such that $K^2=2\chi-6$ or $K^2=2\chi-5$, Gieseker's moduli space $\mathfrak{M}_{K^2,\chi}$ contains surfaces admitting a $\mathbb{Z}_2^2$-action by the results of \cite{Lorenzo2021}. In this note this is extended by the following:

\begin{theorem}\label{TheoremGeographyZ22Actions}
 Let $(K^2, \chi)$ be an admissible pair such that $2\chi-6\leq K^2\leq 8\chi-8$ or $K^2=8\chi$.
 Then $\mathfrak{M}_{K^2, \chi}$ contains surfaces
with a $\mathbb{Z}_2^2$-action.
\end{theorem}

On the other hand, the locus $\overline{\mathfrak{M}}^{Gor}_{K^2, \chi}$ of Gorenstein stable surfaces is known to be open in the moduli space of stable surfaces $\overline{\mathfrak{M}}_{K^2, \chi}$ (cf. \cite[Section 2.B]{FPR17}). Nevertheless, it may happen that 
it is also closed in $\overline{\mathfrak{M}}_{K^2, \chi}$, i.e. it may happen that 
$\overline{\mathfrak{M}}^{Gor}_{K^2, \chi}$ is a union of connected components of $\overline{\mathfrak{M}}_{K^2, \chi}$. In \cite[Example 5.5]{Anthes} it is proved that $\overline{\mathfrak{M}}^{Gor}_{2,4}$ 
is not closed in $\overline{\mathfrak{M}}_{2,4}$.
Considering natural deformations (see \cite[Section 5]{Par1991}) of our constructions of surfaces with a $\mathbb{Z}_2^2$-action in $\mathfrak{M}_{K^2, \chi}$ we are able to show that:
\begin{theorem}\label{GorensteinLocus}
 Let $(K^2, \chi)$ be an admissible pair such that $2\chi-6\leq K^2\leq 8\chi-8$.
 Then $\overline{\mathfrak{M}}^{Gor}_{K^2, \chi}$ is not closed in $\overline{\mathfrak{M}}_{K^2, \chi}$.
\end{theorem}

The note is structured as follows. 
%
In Section \ref{Z22 covers.}
how to construct $\mathbb{Z}_2^2$-covers and to obtain information about them is explained. 
Section \ref{GeographyZ22Section} is devoted to prove 
Theorem \ref{TheoremGeographyZ22Actions} and Theorem \ref{GorensteinLocus}. 
We will construct the surfaces that prove 
Theorem \ref{TheoremGeographyZ22Actions} as $\mathbb{Z}_2^2$-covers of rational surfaces. The surfaces that will be constructed in the region $2\chi-6\leq K^2\leq 4\chi-4$ are either surfaces 
with a genus $2$ fibration or canonical models of surfaces with a genus $2$ fibration.
The surfaces that will be constructed in the region $4\chi-3\leq K^2\leq 8\chi-8$ are surfaces with a genus $3$
fibration. The surfaces that will be constructed on the line $K^2=8\chi$ are products of curves.
The proof of Theorem \ref{GorensteinLocus} consists in finding $\mathbb{Q}$-Gorenstein families $\mathcal{X}\rightarrow T$ of $\mathbb{Z}^2_2$-covers such that:
\begin{enumerate}
 \item[i)] The fiber $\mathcal{X}_t$ is contained in $\overline{\mathfrak{M}}^{Gor}_{K^2,\chi}$ for every $t\in T\setminus \{0\}$.
 \item[ii)] The special fiber $\mathcal{X}_0$ belongs to  $\overline{\mathfrak{M}}_{K^2,\chi}\setminus \overline{\mathfrak{M}}^{Gor}_{K^2,\chi}$.
\end{enumerate}
We will obtain these families considering natural deformations of the $\mathbb{Z}_2^2$-covers that prove 
Theorem \ref{TheoremGeographyZ22Actions}. 
It is worth noticing that the non-Gorenstein
degenerations that will be considered
are normal unless $K^2=2\chi-6$. They have singularities of type
$\frac{1}{4}(1,1)$ (see Remark \ref{1411Singularities} below) when $K^2\neq 2\chi-6$ and they have singularities of type 
$(xy=0)\subset \frac{1}{2}(1,1,1)$ (cf. \cite[Theorem 4.23.iii)]{KSB88})
when $K^2=2\chi-6$.

\section{\texorpdfstring{$\mathbb{Z}_2^2$}{Z22}-covers.}\label{Z22 covers.}
A $\mathbb{Z}_2^2$-cover of a variety $Y$ is a finite map $f:X\rightarrow Y$ together with a faithful action of 
$\mathbb{Z}_2^2$ on $X$ such that $f$ exhibits $Y$ as $X/\mathbb{Z}_2^2$.
The structure theorem for smooth $\mathbb{Z}_2^2$-covers was first 
given by Catanese \cite{Cata1984}.  
According to {\cite[Section 2]{Cata1984}} or {\cite[Theorem 2.1]{Par1991}},
to define a $\mathbb{Z}_2^2$-cover $X\rightarrow Y$ of a
smooth and irreducible projective variety $Y$ with normal $X$, it suffices to consider both:
\begin{enumerate}
 \item[-] Effective divisors $D_1, D_2, D_3$ such that the branch locus $B=D_1+D_2+D_3$ is reduced.
 \item[-] Non-trivial line bundles $L_1, L_2, L_3$ satisfying $2L_1\equiv D_2+D_3, 2L_2\equiv D_1+D_3$ and such that  $L_3\equiv L_1+L_2-D_3$.
\end{enumerate}
The set $\{L_i,D_j\}_{i,j}$ is called the building data of the cover.


\begin{proposition}             
[{{\cite[Section 2]{Cata1984}}} or {{\cite[Proposition 4.2]{Par1991}}}] \label{Z22Formulas}  
Let $Y$ be a smooth surface and $f:X\rightarrow Y$ a smooth $\mathbb{Z}_2^2$-cover with building data 
$\{L_i,D_j\}_{i,j}$. Then:
\begin{equation*}
\begin{split}
2K_X\equiv f^*(2K_Y+D_1+D_2+D_3),\\
K_X^2=(2K_Y+D_1+D_2+D_3)^2,\\
p_g(X)=p_g(Y)+\sum_{i=1}^3h^0(K_Y+L_i),\\
\chi(\mathcal{O}_X)=4\chi(\mathcal{O}_Y)+\frac{1}{2}\sum_{i=1}^3L_i(L_i+K_Y).
\end{split}
\end{equation*}
\end{proposition}

\begin{remark}\label{ChangingDivisorsIsNatural}
 Let $f:X\rightarrow Y$ be a $\mathbb{Z}_2^2$-cover
 of a smooth and irreducible projective variety $Y$ with building data $\{L_i,D_j\}_{i,j}$. Then we can consider a $\mathbb{Z}_2^2$-cover $f':X'\rightarrow Y$ with building data $\{L_i,D'_j\}_{i,j}$ where $D'_j\in |D_j|$ for every $j\in\{1,2,3\}$. Then $f'$ is a natural deformation of $f$ (see \cite[Section 2]{Cata1984} or \cite[Section 5]{Par1991}) and
 we can deform $X$ to $X'$ via a deformation $\phi:\mathcal{X}\rightarrow S$ (cf. \cite[Proposition 5.1]{Par1991}). Moreover, if $S$ is reduced and $\mathcal{X}_s$ is a stable surface for every $s\in S$, then $\phi$ is a $\mathbb{Q}$-Gorenstein deformation by Kollár's numerical criterion for Cartier divisors (see \cite[Theorem 89]{Kollar2016} and \cite[Theorem 27]{KollarMumford}).
\end{remark}

\begin{example}\label{1411Singularities}
Let $r:S\rightarrow X$ be the minimal resolution of a normal surface singularity $(X,p)$
whose exceptional divisor $r^{-1}(p)$ is a $(-4)$-curve $C$. Then $p$ is said to be a $\frac{1}{4}(1,1)$-singularity.
 These singularities are log canonical because $K_S=r^*K_X-\frac{1}{2}C$ (cf. \cite{Alexeev1992}). Remark that
 $\frac{1}{4}(1,1)$-singularities are not Gorenstein because $K_X$ has index $2$.
 
 There is an easy way to obtain 
 $\frac{1}{4}(1,1)$-singularities via $\mathbb{Z}_2^2$-covers (cf. \cite{Catanese1999}). Indeed, 
 let $\pi:X\rightarrow Y$ be a $\mathbb{Z}_2^2$-cover with building data $\{L_i,D_j\}_{i,j}$.
 Let us assume that $D_1, D_2$ and $D_3$ pass through a point $q\in Y$ giving rise to an ordinary triple point
 on the branch locus $B=D_1+D_2+D_3$ of $\pi$. Then $X$ has a $\frac{1}{4}(1,1)$-singularity $p$ over $q$ (see \cite[Proposition 3.3]{Par1991}).
 We can resolve this singularity in a canonical way.
 Let $b:\widetilde{Y}\rightarrow Y$ be the blow-up of $Y$ at $q$ with exceptional divisor $E$. Then there is 
 a $\mathbb{Z}_2^2$-cover $\widetilde{X}\rightarrow\widetilde{Y}$ with branch locus $\widetilde{B}=\widetilde{D_1}+\widetilde{D_2}+
\widetilde{D_3}$ where $\widetilde{D_i}=b^*D_i-E$. The induced map $\widetilde{X}\rightarrow X$ resolves the 
$\frac{1}{4}(1,1)$-singularity. Moreover, using the formulas for $\mathbb{Z}_2^2$-covers it can be proved that 
$\chi(\mathcal{O}_X)=\chi(\mathcal{O}_{\widetilde{X}})$ and $K^2_X=K^2_{\widetilde{X}}+1$.

Now, suppose that we were able to consider a $\mathbb{Z}_2^2$-cover $\pi':X'\rightarrow Y$ branched along $B'=D_1+D_2+D_3'$ where 
$D'_3\in|\mathcal{O}_Y(D_3)|$ does not pass through $q$. Then $X$ would be a degeneration of a surface $X'$ without
$\frac{1}{4}(1,1)$-singularities (see Remark \ref{ChangingDivisorsIsNatural}). Notice that this is consistent with the fact that
$\frac{1}{4}(1,1)$-singularities admit $\mathbb{Q}$-Gorenstein smoothings (cf. \cite[Proposition 2.7]{HackProk}).
\end{example}


\section{Constructions.}\label{GeographyZ22Section}

The aim of this section is to prove Theorem \ref{TheoremGeographyZ22Actions} and Theorem \ref{GorensteinLocus}. Given an admissible pair $(K^2,\chi)$ such that $2\chi-6\leq K^2\leq 8\chi-8$ or $K^2=8\chi$ we are going to construct a surface $S\in\mathfrak{M}_{K^2,\chi}$ with a $\mathbb{Z}_2^2$-action as a $\mathbb{Z}_2^2$-cover of a rational surface. Then we are going to show that, doing slight changes in the branch locus of the covers considered when $2\chi-6\leq K^2\leq 8\chi-8$, we can obtain a $\mathbb{Q}$-Gorenstein degeneration of $S$ with fibers containing non-Gorenstein singularities.


%


\subsection{Genus \texorpdfstring{$2$}{2} fibrations with a \texorpdfstring{$\mathbb{Z}_2^2$}{Z22}-action.}

 In \cite{HorGen2} Horikawa studied genus $2$ fibrations. Amongst the things he showed we have that given a relatively minimal genus $2$
fibration $f:X\rightarrow B$ from a smooth surface $X$, the following formula holds:
\begin{equation*} \label{FormulaGenus2Fibrations}
 K^2_X=2\chi(\mathcal{O}_X)-6\chi(\mathcal{O}_B)+\sum_{p\in B}H(p).
\end{equation*}
The contribution $H(p)$ is a non-negative integer that is bigger than $0$ if and only if 
the fiber over $p\in B$ is not $2$-connected. The possible fibers of a genus $2$ fibration
are well known (e.g. \cite{HorGen2}). 
One of the most relevant not $2$-connected fibers are the so called
fibers of general $(I_1)$-type (see \cite{HorLocal}).
A fiber of general $(I_1)$-type over $p\in B$ consists of $2$ elliptic curves with 
self-intersection $(-1)$ intersecting transversally in one point and $H(p)=1$.

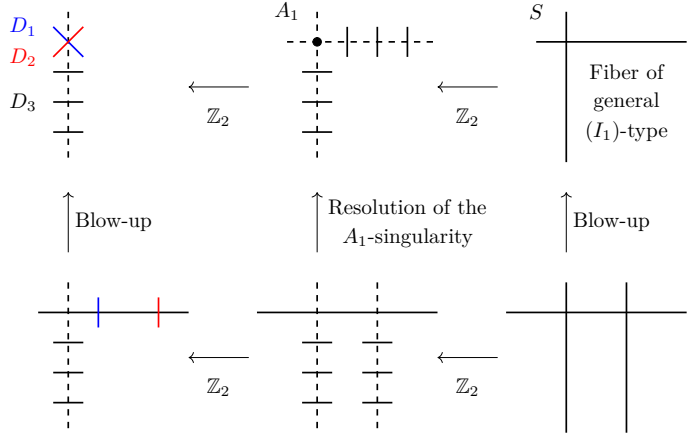
\begin{figure}
\begin{center}
\scalebox{0.8}{
 \begin{tikzpicture}
 \draw[thick, dashed](0,3/2)--(0,-1) ;
 \draw[thick, blue](-1/4,5/4)--(1/4,3/4) node at (-3/4,5/4) {$D_1$}  ;
 \draw[thick, red](-1/4,3/4)--(1/4,5/4) node at (-3/4,3/4) {$D_2$} ;
 \draw[thick](-1/4,1/2)--(1/4,1/2);
 \draw[thick](-1/4,0)--(1/4,0) node at (-3/4,0) {$D_3$};
 \draw[thick](-1/4,-1/2)--(1/4,-1/2);
  \draw[->] (3,1/4)--(2,1/4) node at (5/2,-1/4) {$\mathbb{Z}_2$};
  \draw[->] (0,-5/2)--(0,-3/2)  node at (3/4,-2) {Blow-up};
   \draw[thick, dashed](0,-3)--(0,-11/2) ;
   \draw[thick](-1/2,-7/2)--(2,-7/2) ;
 \draw[thick, blue](1/2,-13/4)--(1/2,-15/4)  ;
 \draw[thick, red](3/2,-13/4)--(3/2,-15/4);
 \draw[thick](-1/4,-8/2)--(1/4,-8/2);
 \draw[thick](-1/4,-9/2)--(1/4,-9/2);
 \draw[thick](-1/4,-10/2)--(1/4,-10/2);
  \draw[->] (3,-17/4)--(2,-17/4) node at (5/2,-19/4) {$\mathbb{Z}_2$};
 \end{tikzpicture}
 \begin{tikzpicture}
 \draw[thick, dashed](0,3/2)--(0,-1) ;
 \draw[thick, dashed](-1/2,1)--(2,1) ;
 \filldraw(0,1) circle (2pt) node at (-1/2,3/2) {$A_1$};
  \draw[thick](1/2,5/4)--(1/2,3/4)  ;
 \draw[thick](3/2,5/4)--(3/2,3/4);
  \draw[thick](1,5/4)--(1,3/4);
 \draw[thick](-1/4,1/2)--(1/4,1/2);
 \draw[thick](-1/4,0)--(1/4,0);
 \draw[thick](-1/4,-1/2)--(1/4,-1/2);
  \draw[->] (3,1/4)--(2,1/4) node at (5/2,-1/4) {$\mathbb{Z}_2$};
  \draw[->] (0,-5/2)--(0,-3/2) node at (3/2,-7/4) {Resolution of the};
   \draw[thick, dashed](0,-6/2)--(0,-11/2) node at (3/2,-9/4) {$A_1$-singularity};
   \draw[thick](-1,-7/2)--(2,-7/2);
 \draw[thick](-1/4,-8/2)--(1/4,-8/2);
 \draw[thick](-1/4,-9/2)--(1/4,-9/2);
 \draw[thick](-1/4,-10/2)--(1/4,-10/2);
    \draw[thick, dashed](1,-6/2)--(1,-11/2) ;
  \draw[thick](3/4,-8/2)--(5/4,-8/2);
 \draw[thick](3/4,-9/2)--(5/4,-9/2);
 \draw[thick](3/4,-10/2)--(5/4,-10/2);
  \draw[->] (3,-17/4)--(2,-17/4) node at (5/2,-19/4) {$\mathbb{Z}_2$};
 \end{tikzpicture}
 \begin{tikzpicture}
 \draw[thick](0,3/2)--(0,-1) node at (1,1/2) {Fiber of};
 \draw[thick](-1/2,1)--(2,1) node at (1,-1/2) {$(I_1)$-type};
  \draw[->] (0,-5/2)--(0,-3/2) node at (3/4,-2) {Blow-up};
   \draw[thick](0,-6/2)--(0,-11/2) node at (1,0) {general};
   \draw[thick](-1,-7/2)--(2,-7/2) node at (-1/2,3/2) {$S$};;
    \draw[thick](1,-6/2)--(1,-11/2) ;
 \end{tikzpicture}
}
 \end{center}
 \caption{Obtaining fibers of general $(I_1)$-type via $\mathbb{Z}_2^2$-covers.} \label{I1typeFigure}
  \end{figure}

Let $(K^2, \chi)$ be an admissible pair such that $2\chi-6\leq K^2\leq 4\chi-4$. In general we are going to obtain a surface 
 $S\in\mathfrak{M}_{K^2, \chi}$ admitting a $\mathbb{Z}_2^2$-action as follows. 
 Let $S\rightarrow \mathbb{F}_e$ be a $\mathbb{Z}_2^2$-cover of the Hirzebruch surface $\mathbb{F}_e$ with negative section 
 $\Delta_0$ of self-intersection $(-e)$ and fiber $F$ whose branch 
locus $B=D_1+D_2+D_3$ is a normal crossing divisor consisting of smooth divisors
\begin{equation*}
\begin{split}
D_1\in |\Delta_0+\alpha F|\\
D_2\in |\Delta_0+\beta F|\\
D_3\in |3\Delta_0+\gamma F|
\end{split}
\end{equation*}
A general fiber of $\mathbb{F}_e$ pulls back to a smooth genus $2$ curve in $S$, but a fiber passing through
a  point of transversal intersection of $D_1$ and $D_2$ becomes a fiber 
of general $(I_1)$-type in $S$ (see Figure \ref{I1typeFigure}). Hence, if we make $D_1$ and $D_2$ intersect transversally 
$\varepsilon := K^2-(2\chi-6)$ times, 
we can obtain a surface $S$ with a genus $2$ fibration over $\mathbb{P}^1$ whose unique not 
 $2$-connected fibers are $\varepsilon$ fibers of general $(I_1)$-type and therefore $K^2_S=2\chi(\mathcal{O}_S)-6+\varepsilon$.
 Proceeding like this we are able to fill in the region $2\chi-6\leq K^2\leq 4\chi-4$ with canonical models
 with the exception of the cases $(K^2, \chi)=(1,2), (1,3)$ and the line $K^2=4\chi-5$. In order to obtain surfaces with a $\mathbb{Z}_2^2$-action in $\mathfrak{M}_{1,2}$ or $\mathfrak{M}_{1,3}$ it suffices to consider $\mathbb{Z}_2^2$-covers of $\mathbb{P}^2$.
 In order to fill in the line $K^2=4\chi-5$ we are going to construct surfaces 
 admitting a $\mathbb{Z}_2^2$-action on the line 
 $K^2=4\chi-4$ and having a $\frac{1}{4}(1,1)$-singularity. Their minimal resolutions will belong to the line $K^2=4\chi-5$ because $\frac{1}{4}(1,1)$-singularities do not affect the holomorphic Euler characteristic  of the minimal resolution but reduce by $1$ the self-intersection of its canonical class (see Example \ref{1411Singularities}). The rest of the subsection is devoted to construct explicitly these surfaces and to obtain non-Gorenstein degenerations of them.
 
  Let us suppose first that $(K^2, \chi)$ is an admissible pair on the line $K^2=2\chi-6$. We can assume that $(K^2, \chi)\neq (2,4)$ because for these invariants Theorem \ref{TheoremGeographyZ22Actions} follows from \cite{Lorenzo2021} and Theorem \ref{GorensteinLocus} follows from \cite[Example 5.5]{Anthes}. We define $\alpha, \beta, \gamma, e$ as follows:
  \begin{enumerate}
   \item[-] If $\chi\equiv 0(2)$ we take $\alpha=0, \beta=2, \gamma=\chi +4, e=2$.
   \item[-] If $\chi\equiv 1(2)$ we take $\alpha=0, \beta=0, \gamma=\chi+1, e=0$.
  \end{enumerate}
  
 In \cite{Lorenzo2021} a 
  surface $S\in\mathfrak{M}_{K^2,\chi}$ with a $\mathbb{Z}_2^2$-action is obtained as a
smooth $\mathbb{Z}_2^2$-cover $\pi:S\rightarrow \mathbb{F}_e$
whose branch 
locus $B=D_1+D_2+D_3$ is a normal crossing divisor consisting of smooth divisors
\begin{equation*}
\begin{split}
D_1\in |\Delta_0+\alpha F|\\
D_2\in |\Delta_0+\beta F|\\
D_3\in |3\Delta_0+\gamma F|
\end{split}
\end{equation*} 
  
We consider now a $\mathbb{Z}_2^2$-cover $\rho:X\rightarrow \mathbb{F}_e$ with branch locus
$B'=D_1+D_2'+D_3$ where:
\begin{enumerate}
 \item[-] If $\chi$ is even $D'_2$ is the union 
of the negative section $\Delta_0=D_1$ and two different general fibers of $\mathbb{F}_2$.
 \item[-] If $\chi$ is odd $D'_2=D_1$.
\end{enumerate}
We notice first that $X$ is not normal because of
\cite[Corollary 3.1]{Par1991}. Now,
the pair $(\mathbb{F}_e, \frac{1}{2}B')$
is log canonical 
because 
$K_{\mathbb{F}_e}+\frac{1}{2}B'$ has index $2$, every prime divisor of $\frac{1}{2}B'$ has multiplicity $\leq 1$ 
and the support of $B'$ is a normal crossing divisor (see \cite[Definition 2.1]{AlexeevPardini2012}). Hence, 
$X$ is semi log canonical by \cite[Proposition 2.5.2]{AlexeevPardini2012}. In addition the canonical divisor of $X$ is ample because 
$2K_X$ is the pullback via $\rho$ of the ample divisor $\Delta_0+(\alpha+\beta+\gamma-2e-4)F$ (cf. \cite[Proposition 2.5.1]{AlexeevPardini2012}).
It follows that $X\in\overline{\mathfrak{M}}_{K^2, \chi}$. Since $X$ is a degeneration of $S$ by Remark \ref{ChangingDivisorsIsNatural},
showing that $X$ is not Gorenstein suffices to prove Theorem \ref{GorensteinLocus} when $K^2=2\chi-6$. Now, 
if we denote by $X^{\nu}\rightarrow X$ the normalization of $X$, we check 
using the normalization algorithm in \cite[Section 3]{Par1991} that $X^{\nu}$ can be obtained as a $\mathbb{Z}_2^2$-cover 
$X^{\nu}\rightarrow\mathbb{F}_e$
whose branch 
locus $C=C_1+C_2+C_3$ consists of
\begin{equation*}
\begin{split}
C_1=0\\
C_2=D'_2-D_1\\
C_3=D_1+D_3
\end{split}
\end{equation*}
Note that 
$X^{\nu}$ consists of two disjoint copies of a surface with $A_1$-singularities over the singularities of $C$ (cf. \cite[Proposition 3.1]{Par1991}).
The surface $X$ is obtained by gluing the two components of $X^{\nu}$ along the curve over $D_1$.
The singularities on $X$ over the points of intersection of $D_1$ and $D_3$ are of the so called type $(xy=0)\subset \frac{1}{2}(1,1,1)$ (cf. \cite[Theorem 4.23.iii)]{KSB88})
and they have index $2$.
Hence 
$X\notin\overline{\mathfrak{M}}^{Gor}_{K^2, \chi} $ and our claim follows.

  When $(K^2, \chi)=(1,2)$ we can obtain a surface 
  $S\in\mathfrak{M}_{K^2,\chi}$ with a $\mathbb{Z}_2^2$-action as a $\mathbb{Z}_2^2$-cover 
$\pi:S\rightarrow \mathbb{P}^2$ whose  branch locus
$B = D_1 + D_2 + D_3$ is a normal crossing divisor consisting of a line $D_1$, and smooth cubics $D_2,D_3$.
Let us consider a $\mathbb{Z}_2^2$-cover $X\rightarrow \mathbb{P}^2$ with branch locus 
$B' = D'_1 + D_2 + D_3$ where $D'_1$ is a line passing through one and only one point of intersection $p$ of $D_2$ and $D_3$
and intersects them transversally everywhere. Then $X$ has a $\frac{1}{4}(1,1)$-singularity over $p$, $K_X$ is ample, 
$\chi(\mathcal{O}_X)=2$ and $K^2_X=1$. Moreover, $X$ is a degeneration of $S$ by Remark \ref{ChangingDivisorsIsNatural}
and therefore this proves Theorem \ref{GorensteinLocus} when $(K^2, \chi)=(1,2)$.

When $(K^2, \chi)=(1,3)$ we can obtain a surface 
  $S\in\mathfrak{M}_{K^2,\chi}$ with a $\mathbb{Z}_2^2$-action as a $\mathbb{Z}_2^2$-cover  $\pi:S\rightarrow \mathbb{P}^2$ whose  branch locus
$B = D_1 + D_2 + D_3$ is a normal crossing divisor consisting of two
different lines $D_1 , D_2$ and a quintic $D_3$.
Let us consider a $\mathbb{Z}_2^2$-cover $X\rightarrow \mathbb{P}^2$ with branch locus 
$B' = D'_1 + D_2 + D_3$ where $D'_1$ is a line passing through one and only one point of intersection $p$ of $D_2$ and $D_3$
and intersects them transversally everywhere. Then $X$ has a $\frac{1}{4}(1,1)$-singularity over $p$, $K_X$ is ample, 
$\chi(\mathcal{O}_X)=3$ and $K^2_X=1$. Moreover, $X$ is a degeneration of $S$ by Remark \ref{ChangingDivisorsIsNatural}
and therefore this proves Theorem \ref{GorensteinLocus} when $(K^2, \chi)=(1,3)$.

  Let us fix now an admissible pair $(K^2, \chi)\neq (1,2), (1,3)$ such that $2\chi-5\leq K^2 \leq  4\chi-6$.
  We define $\alpha, \beta, \gamma, e$ as follows:
  \begin{enumerate}
   \item[-] If $K^2\equiv 0(4)$ we take $\alpha=0, \beta=K^2-2\chi+6, \gamma=2\chi - 2-\frac{1}{2}K^2, e=0$.
   \item[-] If $K^2\equiv 2(4)$ we take $\alpha=1, \beta=K^2-2\chi+5, \gamma=2\chi - 2-\frac{1}{2}K^2, e=0$.
   \item[-] If $K^2\equiv 3(4)$ we take $\alpha=0, \beta=K^2-2\chi+7, \gamma=2\chi - \frac{1}{2}-\frac{1}{2}K^2, e=1$.
   \item[-] If $K^2\equiv 1(4)$ we take $\alpha=1, \beta=K^2-2\chi+6, \gamma=2\chi - \frac{1}{2}-\frac{1}{2}K^2, e=1$.
  \end{enumerate}
  
  We consider a smooth $\mathbb{Z}_2^2$-cover $\pi:S\rightarrow \mathbb{F}_e$ whose branch 
locus $B=D_1+D_2+D_3$ is a normal crossing divisor consisting of smooth divisors
\begin{equation*}
\begin{split}
D_1\in |\Delta_0+\alpha F|\\
D_2\in |\Delta_0+\beta F|\\
D_3\in |3\Delta_0+\gamma F|
\end{split}
\end{equation*}
According to the formulas for $\mathbb{Z}_2^2$-covers 
 (see Proposition \ref{Z22Formulas})
 we have:
 \begin{equation*}
 \begin{split}
   \chi(\mathcal{O}_S)=\frac{1}{2}(\alpha+\beta)+\gamma-2e-1=\chi,\\
   2K_S=\pi^*\left(\Delta_0+\frac{1}{2}(K^2+e)F\right),\\
   K^2_S=2\alpha+2\beta+2\gamma-5e-8=K^2.
 \end{split}
 \end{equation*}
In particular $K_S$ is ample and Theorem \ref{TheoremGeographyZ22Actions} is proved when
 $2\chi-6\leq K^2\leq 4\chi-6$. Moreover, 
\begin{equation*}
 \begin{split}
D_1D_2=K^2-2\chi+6>0,\\
D_1D_3=3\alpha+\gamma-3e>0,\\
D_2D_3=D_1D_3+3(\beta-\alpha)>0,\\
h^0(\mathcal{O}_{\mathbb{F}_0}(D_3))=8\chi-4-2K^2\geq 8.
 \end{split}
\end{equation*}
It follows that we can find a smooth $D'_3\in|D_3|$ passing through a unique point of intersection $p$ of $D_1$ and $D_2$
and intersecting these divisors transversally everywhere. If $X\rightarrow \mathbb{F}_e$ is a $\mathbb{Z}_2^2$-cover 
 with branch locus $B'=D_1+D_2+D'_3$, then $K_X$ is ample, $\chi(\mathcal{O}_X)=\chi$ and $K^2_X=K^2$. Moreover, 
 the only singularity of $X$ is a $\frac{1}{4}(1,1)$-singularity over $p$ (see Example \ref{1411Singularities}).
 Since $X$ is a degeneration of $S$ by 
 Remark \ref{ChangingDivisorsIsNatural}, this proves Theorem \ref{GorensteinLocus} when
 $2\chi-5\leq K^2\leq 4\chi-6$ and $(K^2,\chi)\neq(1,2), (1,3)$.

We conclude that in order to prove Theorem \ref{TheoremGeographyZ22Actions} and Theorem \ref{GorensteinLocus} for pairs in the region 
$2\chi-6\leq K^2\leq 4\chi-4$
the only cases left to study are $K^2=4\chi-5$
and $K^2=4\chi-4$.

Let us consider a $\mathbb{Z}_2^2$-cover $\pi:T\rightarrow \mathbb{F}_0$ 
of $\mathbb{F}_0$ whose branch 
locus $B=D_1+D_2+D_3$ consists of divisors
\begin{equation*}
\begin{split}
D_1\in |\Delta_0+2 F|\\
D_2\in |\Delta_0+2\chi F|\\
D_3\in |3\Delta_0|
\end{split}
\end{equation*}
We take $D_1$ and $D_2$ smooth and irreducible intersecting transversally in $D_1\cdot D_2= 2\chi+2$ points and 
$D_3$ consisting of $3$ different fibers
$\Delta_1,\Delta_2,\Delta_3$ such that only $\Delta_1$ passes through a point of intersection 
$p$ of $D_1$ and $D_2$. The surface $T$ has a $\frac{1}{4}(1,1)$-singularity over $p$, 
$\chi(\mathcal{O}_T)=\chi$ and $K^2_T=K^2+1$. If we denote by $q:\widetilde{\mathbb{F}}_0\rightarrow \mathbb{F}_0$ the blow-up of $\mathbb{F}_0$ at $p$ with exceptional divisor $E$ and $S\rightarrow T$ is the minimal resolution  of $T$, then $\pi$ induces a $\mathbb{Z}_2^2$-cover $\widetilde{\pi}:S\rightarrow \widetilde{\mathbb{F}}_0$
with branch locus $\widetilde{B}=\widetilde{D}_1+\widetilde{D}_2+\widetilde{D}_3$ where 
$\widetilde{D}_i=q^*D_i-E$ for each $i\in\{1,2,3\}$ and 
$K^2_S=K^2, \chi(\mathcal{O}_{S})=\chi$ (see Example \ref{1411Singularities}). Moreover, $K_S$ is ample because $2K_S$ is the pullback via $\widetilde{\pi}$ of $D:=q^*(\Delta_0+(2\chi-2)F)-E$ and this divisor can be checked to be ample using Nakai-Moishezon criterion. 
Indeed, let us suppose that there exists an irreducible curve $C\in |q^*(a\Delta_0+bF)-cE|$ such that $CD<0$ for some non-negative integers $a,b,c$. This implies that $q(C)\in |a\Delta_0+bF|$ is an irreducible curve of $\mathbb{F}_{0}$ with a point of multiplicity $c>(2\chi-2)a+b$, which is clearly impossible.
Therefore these examples prove Theorem \ref{TheoremGeographyZ22Actions} when $K^2=4\chi-5$.
We define now $D'_3:=\Delta_1+\Delta_2+\Delta'_3$ where $\Delta'_3\in|\Delta_0|$ passes
through a point of intersection $p'$ of $D_1$ and $D_2$ different from 
$p$ and $\widetilde{D}'_3=q^*D'_3-E\in |\widetilde{D}_3|$. If $X\rightarrow \widetilde{\mathbb{F}}_0$ is the $\mathbb{Z}_2^2$-cover of $\widetilde{\mathbb{F}}_0$ with branch locus 
$\widetilde{B}'=\widetilde{D}_1+\widetilde{D}_2+\widetilde{D}'_3$
then $K_X$ is ample, $\chi(\mathcal{O}_X)=\chi$ and $K^2_X=K^2$. Moreover, 
 the only singularity of $X$ is a $\frac{1}{4}(1,1)$-singularity coming from $p'$ (see Example \ref{1411Singularities}).
 Since $X$ is a degeneration of $S$ by 
 Remark \ref{ChangingDivisorsIsNatural}, this proves Theorem \ref{GorensteinLocus} when
 $K^2=4\chi-5$.
 
 If we want to obtain a surface that admits a $\mathbb{Z}_2^2$-action in $\mathfrak{M}_{K^2, \chi}$ when $K^2=4\chi-4$ it suffices to consider a smooth $\mathbb{Z}_2^2$-cover $\pi:S\rightarrow \mathbb{F}_0$ 
 whose branch locus $B=D_1+D_2+D_3$ is a normal crossing divisor consisting of smooth divisors
\begin{equation*}
\begin{split}
D_1\in |\Delta_0+2 F|\\
D_2\in |\Delta_0+2\chi F|\\
D_3\in |3\Delta_0|
\end{split}
\end{equation*}
Indeed, according to the formulas for $\mathbb{Z}_2^2$-covers 
 (see Proposition \ref{Z22Formulas})
 we have $\chi(\mathcal{O}_S)=\chi$ and $K^2_S=K^2$.
 Moreover, $K_S$ is ample because the divisor $2K_S$ is the pullback via $\pi$ of the ample divisor $\Delta_0+(2\chi-2)F$. Let $D'_3\in|D_3|$ consist of $3$ different fibers such that one and only one of them passes through a point of intersection $p$ of $D_1$ and $D_2$. If $X\rightarrow \mathbb{F}_0$ is a $\mathbb{Z}_2^2$-cover 
 with branch locus $B'=D_1+D_2+D'_3$, then $K_X$ is ample, $\chi(\mathcal{O}_X)=\chi$ and $K^2_X=K^2$. Moreover, 
 the only singularity of $X$ is a $\frac{1}{4}(1,1)$-singularity over $p$ (see Example \ref{1411Singularities}).
 Since $X$ is a degeneration of $S$ by 
 Remark \ref{ChangingDivisorsIsNatural}, this proves Theorem \ref{GorensteinLocus} when
 $K^2=4\chi-4$. 

%
%

\subsection{Genus \texorpdfstring{$3$}{3} fibrations with a \texorpdfstring{$\mathbb{Z}_2^2$}{Z22}-action.}
 \label{ConstructionsGenus3WithZ22Actions}
 In order to fill in the region $4\chi-3\leq K^2\leq 8\chi-8$ with surfaces admitting a
 $\mathbb{Z}_2^2$-action we are going to consider surfaces $T\in \mathfrak{M}_{K^2, \chi}$
 that have a genus $3$ fibration. If we want the quotient $T/\mathbb{Z}_2^2$ to be a Hirzebruch surface, the Riemann-Hurwitz formula implies that it can be done as follows. We consider a $\mathbb{Z}_2^2$-cover
 $T\rightarrow \mathbb{F}_e$  of the Hirzebruch surface $\mathbb{F}_e$ with negative section $\Delta_0$ of self-intersection $(-e)$
and fiber $F$ whose branch 
locus $B=D_1+D_2+D_3$ is a normal crossing divisor consisting of smooth divisors
\begin{equation*}
\begin{split}
D_1\in |\alpha F|\\
D_2\in |2\Delta_0+\beta F|\\
D_3\in |4\Delta_0+\gamma F|
\end{split}
\end{equation*}
Nevertheless, in this way we can only obtain surfaces with $K^2\equiv 0(4)$.
To avoid this problem, let $(K^2, \chi)$ be an admissible pair in the region $4\chi-3\leq K^2\leq 8\chi-8$ and define $\varepsilon\in\{0,1,2,3\}$ such that 
$K^2+\varepsilon\equiv 0 (4)$. Then we can construct a $\mathbb{Z}_2^2$-cover 
$T\rightarrow \mathbb{F}_e$ such that $\chi(\mathcal{O}_T)=\chi$ and $K^2_T=K^2+\varepsilon$. If we make the divisor $D_1$ of the cover above pass through $\varepsilon$ points of intersection of $D_2$ and $D_3$, we can produce 
$\varepsilon$ singularities of type $\frac{1}{4}(1,1)$ on $T$  (see Example \ref{1411Singularities}). Considering the minimal resolution of $T$ we will be able to obtain a surface with the invariants that we are interested in because each $\frac{1}{4}(1,1)$-singularity does not affect the holomorphic Euler characteristic  of the minimal resolution but reduces by $1$ the self-intersection of its canonical class. The rest of the subsection is devoted to construct explicitly these surfaces and to obtain non-Gorenstein degenerations of them.

We fix an admissible pair $(K^2, \chi)$ such that $4\chi-3\leq K^2\leq 8\chi-8$ and
$\varepsilon\in\{0,1,2,3\}$ such that $K^2+\varepsilon\equiv 0 (4)$. In addition, we choose $\alpha, \beta, \gamma$ as follows:
\begin{enumerate}
 \item[-] If $K^2+\varepsilon\equiv 0 (8)$ we take $\alpha=\frac{1}{2}(K^2+\varepsilon)-2\chi+4, \beta=2\chi-\frac{1}{4}(K^2+\varepsilon)$ and $\gamma=0$.
 \item[-] If $K^2+\varepsilon\equiv 4 (8)$ we take $\alpha=\frac{1}{2}(K^2+\varepsilon)-2\chi+5, \beta=2\chi-\frac{1}{4}(K^2+\varepsilon)-2$ and $\gamma=1$.
\end{enumerate}
We consider a $\mathbb{Z}_2^2$-cover $\pi:T\rightarrow \mathbb{F}_0$ 
whose branch locus $B=D_1+D_2+D_3$ consists of divisors
\begin{equation*}
\begin{split}
D_1\in |\alpha F|\\
D_2\in |2\Delta_0+\beta F|\\
D_3\in |4\Delta_0+\gamma F|
\end{split}
\end{equation*}
Note that by the choices of $\beta$ and $\gamma$ we can take $D_2$ smooth and intersecting transversally $D_3$ in $D_2D_3\geq 6$ points in general position. In addition, we can make $D_1$ consist of $\alpha\geq 4$ fibers
 $F_1,\ldots, F_{\alpha}$ such that each $F_i, i\in\{1,\ldots,\varepsilon\}$ passes through a point  of intersection $p_i$ of $D_2$ and $D_3$. 
According to the formulas for $\mathbb{Z}_2^2$-covers  (see Proposition \ref{Z22Formulas}) we have:
 \begin{equation*}
  \begin{split}
   \chi(\mathcal{O}_T)=\frac{1}{2}\alpha+\beta+\frac{3}{2}\gamma-2=\chi,\\
   2K_T=\pi^*\left(2\Delta_0+\frac{1}{4}(K^2+\varepsilon)F\right),\\
   K^2_T=4(\alpha+\beta+\gamma)-16=K^2+\varepsilon.
  \end{split}
 \end{equation*}
In addition $T$ has $\varepsilon$ singularities of type $\frac{1}{4}(1,1)$ over $p_1,\ldots, p_{\varepsilon}$. Let $q:\widetilde{\mathbb{F}}_0\rightarrow \mathbb{F}_0$ be the blow up of $\mathbb{F}_0$ at $p_1,\ldots, p_{\varepsilon}$ with exceptional divisors $E_1,\ldots,E_{\varepsilon}$ respectively. Denote by
$\widetilde{\pi}:S \rightarrow \widetilde{\mathbb{F}}_0$ the $\mathbb{Z}_2^2$-cover
with branch locus $\widetilde{B}=\widetilde{D}_1+\widetilde{D}_2+\widetilde{D}_3$ where 
$\widetilde{D}_j=q^*D_j-\sum_{i=1}^{\varepsilon}E_i$ for each $j\in\{1,2,3\}$. Then the induced map $S\rightarrow T$ is the minimal resolution of $T$ and 
$\chi(\mathcal{O}_S)=\chi, K^2_S=K^2$. Moreover, $K_S$ is ample because $2K_S$ is the pullback via $\widetilde{\pi}$ of 
\begin{equation*}
 q^*\left(2\Delta_0+\frac{1}{4}(K^2+\varepsilon)F\right)-\sum_{i=1}^{\varepsilon}E_i
\end{equation*}
and this divisor is ample by the Nakai-Moishezon criterion. These examples prove Theorem \ref{TheoremGeographyZ22Actions} when $4\chi-3\leq K^2\leq 8\chi-8$. 

We define now $D'_1:=\sum_{i=1}^{\alpha-1}F_i+F'_\alpha$ where $F'_\alpha\in|F|$ passes
through a point of intersection $p_{\varepsilon+1}$ of $D_2$ and $D_3$ different from 
$p_1,\ldots, p_{\varepsilon}$ and $\widetilde{D}'_1=q^*D'_1-\sum_{i=1}^{\varepsilon}E_i\in |\widetilde{D}_1|$. If $X\rightarrow \widetilde{\mathbb{F}}_0$ is the $\mathbb{Z}_2^2$-cover of $\widetilde{\mathbb{F}}_0$ with branch locus 
$\widetilde{B}'=\widetilde{D}'_1+\widetilde{D}_2+\widetilde{D}_3$
then $K_X$ is ample, $\chi(\mathcal{O}_X)=\chi$ and $K^2_X=K^2$. Moreover, 
 the only singularity of $X$ is a $\frac{1}{4}(1,1)$-singularity coming from $p_{\varepsilon+1}$ (see Example \ref{1411Singularities}).
 Since $X$ is a degeneration of $S$ by 
 Remark \ref{ChangingDivisorsIsNatural}, this proves Theorem \ref{GorensteinLocus} when
$4\chi-3\leq K^2\leq 8\chi-8$.
 
 \begin{remark}
The general fiber of the genus $3$ fibration that $|F|$ induces on the surfaces $S\in\mathfrak{M}_{K^2, \chi}$ with a $\mathbb{Z}_2^2$-action that we have just constructed is hyperelliptic.
 In addition, \cite[Remark 3.16]{MendesPardini2021} can be used to show that the canonical map of our examples in the region $8\chi-15\leq K^2\leq 8\chi-8$ is composite with a pencil.
\end{remark}

\begin{remark}
  Using the formulas for $\mathbb{Z}_2^2$-covers we can check that all the examples of surfaces $S\in\mathfrak{M}_{K^2, \chi}$ with a $\mathbb{Z}_2^2$-action that we have constructed so far have irregularity $q(S)=0$.
\end{remark}

\subsection{Products of curves with a \texorpdfstring{$\mathbb{Z}_2^2$}{Z22}-action.}\label{ConstructionsProductsWithZ22Actions}
Let $(K^2, \chi)$ be an admissible pair on the line $K^2=8\chi$. Then we can consider the product of a curve of genus $2$
and a hyperelliptic curve of genus $\chi+1$. If we denote by $\Delta_0$ and $F$ the two classes of fibers of 
$\mathbb{P}^1\times \mathbb{P}^1$,
it can be obtained as a smooth $\mathbb{Z}_2^2$-cover $\pi:S\rightarrow \mathbb{P}^1
 \times \mathbb{P}^1$ whose branch locus $B=D_1+D_2+D_3$ consists of smooth divisors
\begin{equation*}
\begin{split}
D_1\in |6\Delta_0|\\
D_2\in |(2\chi+4) F|\\
D_3=0
\end{split}
\end{equation*}
According to Proposition \ref{Z22Formulas} we have $\chi(\mathcal{O}_S)=\chi$ and $K^2_S=K^2$. Moreover, the divisor $K_S$ is ample
because $2K_S$ is the pullback via $\pi$ of the ample divisor $2\Delta_0+2\chi F$. Hence, these examples prove 
Theorem \ref{TheoremGeographyZ22Actions} when $K^2=8\chi$.
\medbreak

\noindent \begin{acknowledgements}
The author is deeply indebted to his supervisor Margarida Mendes Lopes for all her help.
\end{acknowledgements}

\bibliographystyle{plain}      
\bibliography{GeographyZ22}

\newpage

\noindent Vicente Lorenzo \footnote{The author is a collaborator of the Department of Mathematics and
Center for Mathematical Analysis, Geometry and Dynamical Systems of Instituto Superior T\'{e}cnico,
Universidade de Lisboa and was supported by Fundac\~{a}o para a Ci\^{e}ncia e a Tecnologia (FCT), Portugal through
the program Lisbon Mathematics PhD (LisMath), scholarship  FCT - PD/BD/128421/2017 and
projects UID/MAT/04459/2019 and UIDB/04459/2020.}\\
Center for Mathematical Analysis, Geometry and Dynamical Systems\\
Departamento de Matem\'{a}tica\\
Instituto Superior T\'{e}cnico\\
Universidade de Lisboa\\
Av. Rovisco Pais\\
1049-001 Lisboa\\
Portugal\\
\textit{E-mail address: }{vicente.lorenzo@tecnico.ulisboa.pt}\\
https://orcid.org/0000-0003-2077-6095

\end{document}